\documentclass[12pt]{article}

\usepackage{amsmath}
\usepackage{amsthm}
\usepackage{amsfonts}

\newcommand{\Rset}{\mathbb{R}}
\newcommand{\Cset}{\mathbb{C}}

\newtheorem{prop}{PROPOSITION}
\newtheorem{cor}{COROLLARY}
\theoremstyle{remark}
\newtheorem{rem}{REMARK}
\begin{document}
\noindent
\textbf{ICM 2002 Satellite Conference}\\
ABSTRACT AND APPLIED ANALYSIS \\
Hanoi, August 13-17, 2002.

\medskip
\title{The random integral representation hypothesis revisited : new
classes of s-selfdecomposable laws.}

\date{appeared in cf. ABSTRACT AND APPLIED ANALYSIS; Proc. Int. Conf. Hanoi, 13-17 August
2002, WORLD SCIENTIFIC, Hongkong 2004, pp. 495-514.}

\author{Zbigniew J. Jurek\footnote{This work
was completed while the author was a Visiting Professor at Wayne State University, Detroit, USA .} ,
Wroc\l aw, Poland.}
\maketitle

\begin{quote}
\noindent {\footnotesize \textbf{ABSTRACT.} For $\,0<\alpha\le \infty$, new subclasses
$\,\mathcal{U}^{<\alpha>}$ of the class $\,\mathcal{U}$, of
s-selfdecomposable probability measures, are studied. They are
described by random integrals, by their
characteristic functions and their L\'evy spectral measures. Also their
relations with the classical L\'evy class $L$ of selfdecomposable
distributions are investigated.

}
\end{quote}

Limit distribution theory belongs to the core of probability and
mathematical statistics. Often limit laws are described by
\emph{analytical tools} such as Fourier or Laplace transforms,
but a more stochastic approach (e.g., like stochastic integration,
stopping times, random functionals etc.), seems more natural for
probability questions. Some illustrations of this paradigm are given
in the last paragraph of
this note. In a similar spirit, in Jurek (1985) on page 607 (and
later repeated in Jurek (1988) on page 474), the following hypothesis was
formulated:

\medskip
\emph{Each class of limit distributions, derived from sequences of
independent random variables, is the image of some subset of ID
\emph{(the infinitely divisible probability measures)} by some mapping
defined as a random integral.}
\medskip

Random integral representations, when they can be established, would
provide  descriptions of limiting laws via \emph{stochastic
methods}, i.e., as the probability distributions of the random
integrals of form
\begin{equation}
\int_{(a,b]}h(t)\,dY(r(t))=\int_{(r(a),r(b)]}h(r^{-1}(s))\,dY(s),
\end{equation}
where $h$ and $r$ are deterministic functions, $h:(a,b] \to
\Rset$, $r:(a,b] \to (0,\infty)$ and  $Y(s), 0\le s< \infty$, is a
stochastic process with independent and stationary increments and
cadlag (right continuous with left hand limits) paths; in short,
we refer to $Y$ as \emph{a L\'evy process}. In this note we
provide new examples of classes of limit distributions for which
the above hypothesis holds true. The main results here are
Propositions 3, 4 and 5, and Corollaries 5, 6 and 7.

\medskip
\textbf{1. Introduction and notation.}

Let $E$ denotes a real separable Banach space, $E^{\prime}$
its conjugate space, $<
\cdot,\cdot >$ the usual pairing between $E$ and $E^{\prime}$, and
$||.||$ the norm on E. The $\sigma$-field of all Borel subsets of
$E$ is denoted by $\mathcal{B}$, while $\mathcal{B}_{0}$ denotes
Borel subsets of $E \setminus \{0\}$. By $\mathcal{P}(E)$ we
denote the (topological) semigroup of all Borel probability
measures on $E$, with convolution ``$\ast$" and the weak topology,
in which convergence is denoted by ``$\Rightarrow$". Similarly, by
$ID(E)$ we denote the
topological convolution semigroup of all \emph{infinitely
divisible} probability measures, i.e.,
\[
\mu \in ID(E) \ \ \mbox{iff} \ \ \forall(\mbox{natural} \ \ k \ge
2) \ \ \exists( \mu_{k} \in \mathcal{P}(E)) \ \ \mu=\mu_{k}^{\ast
k}.
\]
Recall also here that $ID(E)$ is a closed topological subsemigroup
of $\mathcal{P}(E)$. Finally on a Banach space $E$ we define the
transforms $T_{r}$, for $r>0$, as follows: $T_{r}x:=rx, \ x \in
E$, and define $\mathcal{L}(\xi)$ as the probability distribution
of an E-valued random variable $\xi$.

A probability measure $\mu \in \mathcal{P}(E)$ is said to be
\emph{s-selfdecomposable} on $E$, and we will write $\mu \in
\mathcal{U}(E)$, if there exists a sequence $\rho_n \in ID(E)$
such that
\begin{equation}
\nu_n:=T_{\frac{1}{n}}(\rho_1 \ast \rho_2 \ast ... \ast
\rho_n)^{\ast 1/n} \Rightarrow \mu, \ \ \mbox{as} \ \ n\to \infty.
\end{equation}
Since we begin with infinitely divisible measures $\rho_n$ we do
not include the shifts $\delta_{x_n}$ in (1), and do not assume that the
triangle system $\{T_{\frac{1}{n}}\,\rho_j^{\ast 1/n}: 1\le j\le n;
n\ge 1\}$ is uniformly infinitesimal, as is usually done in the
general limiting distribution theory.
\newline Also let us note that our
definition (2) is, in fact, the result of Theorem 2.5 in Jurek (1985).
There s-selfdecomposability was defined in many different but
equivalent forms. Finally, s-selfdecomposable distributions
appeared in the context of an approximation of processes by their
discretization; cf. Jacod, Jakubowski and M\'emin (2001).

Originally the s-selfdecomposable distributions were introduced as
limit distributions for sums of shrunken random variables in Jurek
(1981). The \emph{'s'} stands here for \emph{shrinking operation}
defined as follows:
\[
U_r(x):=max(||x||-r,0)\,x/||x||, \ \ \ \ \mbox{for} \ \ \ r>0 \ \
\ \mbox{and} \ \ \ x\in E\setminus\{0\}.
\]
Also see the announcement in Jurek (1977). On the real line
similar distributions, but not related to s-operation, were
studied in O'Connor (1979).

In the present paper we will repeat the scheme (2) successively
and will assume that $\rho_k$ are chosen from a previously
obtained class of limit laws. Such an approach, for another scheme
of limiting procedure was introduced by K. Urbanik (1973) and then
continued by K. Sato, A. Kumar and B. M. Schreiber, N. Thu, with
the most general setting, up to now, described in Jurek (1983),
where there is also a list of related references.

\medskip
For easy reference we collect below some of the known
characterizations of the class $\mathcal{U}(E)$ of
s-selfdecomposable probability measures and indicate only the main
steps in the corresponding proofs.
\begin{prop}
The following statements are equivalent:
\begin{itemize}
\item[\emph{(i)}] $\mu \in \mathcal{U}(E)$.

\item[\emph{(ii)}] $\forall(0<c<1) \ \exists(\mu_c \in ID(E)) \ \mu =
T_c\mu^{\ast c} \ast \mu_c $.

\item[\emph{(iii)}] there exists a unique L\'evy process Y such
that  $\mu = \mathcal{L}(\int_{(0,1)}t\,dY(t))$ .
\end{itemize}
\end{prop}
\emph{Sketch of proofs.} Characterizations (i) and (ii) are
equivalent by Theorem 2.5 and Corollary 2.3 in Jurek (1985).
Equivalence of (ii) and (iii) follows from Theorem 1.1 and Theorem
1.2(a) in Jurek (1988), where one needs to take the constat
$\beta=1$ and the linear operator $Q=I$.
\medskip

For our purposes we define \emph{random integrals} by the formal
formula of integration by parts:
\[
\int_{(a,b]}h(t)\,dY(r(t)):= h(b)Y(r(b))-h(a)Y(r(a)) -
\int_{(a,b]}Y(r(t))\,dh(t),
\]
where the later integral is defined as a limit of the appropriate
Rieman-Stieltjes partial sums. This "limited" approach to
integration is sufficient for our purposes; cf. Jurek and Vervaat
(1983) or Jurek and Mason (1993), Section 3.6. On the other hand,
since L\'evy processes are semi-martingales, the integrals (1) or
the above, can be defined as the stochastic integrals as well.
\begin{cor}
The class $\mathcal{U}$ of s-selfdecomposable probability measures
is closed topological convolution subsemigroup of $ID$. Moreover,
it also is closed under the convolution powers (i.e, for $t>0$ and
$\mu$ we have that $ \mu\in \mathcal {U}$ if and only if
$\mu^{\ast t}\in \mathcal{U}$) and the dilations $T_d$, for $d \in
\Rset$ ( i.e., $\mu \in \mathcal{U}$ if and only if $T_d \mu \in
\mathcal{U}$).
\end{cor}
\emph{Proof}. Both algebraic properties follow from (ii) in
Proposition 1 and the following identities $(T_d(\nu \ast
\rho))^{\ast t} = T_d \nu^{\ast t} \ast T_d \rho^{\ast t}$, for
$t>0$, $d \in \Rset$, and $\nu, \rho \in ID$. To show that
$\mathcal{U}$ is closed in weak convergence topology we use again
the factorization (ii) together with Theorem 1.7.1 in Jurek an
Mason (1993) or cf. Chapter 2 in Parthasarathy (1967).

\medskip
In view of the property (iii), in Proposition 1, we define  the
following integral mapping

\begin{equation}
\mathcal{J}: ID(E) \to \mathcal{U}(E) \ \ \mbox{given by} \ \
\mathcal{J}(\rho):= \mathcal{L}(\int_{(0,1)}s\, dY_\rho(s)),
\end{equation}
where $Y_{\rho}(\cdot)$ is a L\'evy process (i.e., a process with
independent and stationary increments, starting from zero and with
cadlag paths) such that $\mathcal{L}(Y_{\rho}(1))= \rho$. We refer
to $Y(\cdot)$ as the \textbf{b}ackground \textbf{d}riving
\textbf{L}\'evy \textbf{p}rocess (in short, the BDLP) for the
s-selfdecomposable measure $\mathcal{J}(\rho)$.
\begin{rem}
The random integral mapping $\mathcal{J}$ is an isomorphism
between the closed topological semigroups $ID(E)$ and
$\mathcal{U}(E)$; cf. Jurek (1985), Theorem 2.6.
\end{rem}
Finally, let
\[
\hat{\mu}(y):= \int_E e^{i<y,x>} \mu(dx),\ \ y \in E^{\prime},
\]
be \emph{the characteristic function} (the Fourier transform) of a
measure $\mu$. Then for random integrals (1) we infer that
\begin{equation}
\Big(\mathcal{L}(\int_{(a,b]}h(t)dY_{\rho}(r(t)))\Big)^{\widehat{
}}(y)= \exp \int_{(a,b]}\log \hat{\rho}(h(t)y)dr(t),
\end{equation}
when $h$ is a deterministic function, $r$ is an increasing (or
monotone) time change in $(0,\infty)$ and $Y_{\rho}(.)$ a L\'evy
process; cf. Lemma 1.2 in Jurek and Vervaat (1983) or Lemma 1.1 in
Jurek (1985) or simply approximate the right-hand integral by
Rieman-Stieltjes partial sums.

Our results are given in the generality of a Banach space $E$, however,
below in many formulas we will skip the dependence on $E$.
\medskip

\textbf{2. m-times s-selfdecomposable probability measures.}

 Let us put $\mathcal{U}^{<1>}:= \mathcal{U}(E)$ and for $m \ge 2$, let
$\mathcal{U}^{<m>}$ denotes the class of limiting measures in (2), when
$\rho_k \in \mathcal{U}^{<m-1>}$, for $k=1,2,...$ . As a
convention we assume that $\mathcal{U}^{<0>}:= ID$. Our first
characterization is proved along the lines of the proofs of
Theorem 1.1 and 1.2 in Jurek (1988), however one needs not to
confuse the classes $\mathcal{U}_{\beta}$ introduced there, with
those of $\mathcal{U}^{<m>}$ investigated here. Needed changes in
arguments are explained as they are deemed.

\begin{prop}
For $m=1,2,... $, the following are equivalent descriptions of
m-times s-selfdecomposable probability measures:

\begin{itemize}
\item[\emph{(i)}] $\mu \in \mathcal{U}^{<m>}$.

\item[\emph{(ii)}] $\forall(0<c<1) \  \exists(\mu_c \in \mathcal{U}^{<m-1>}) \ \mu =
T_c\mu^{\ast c} \ast \mu_c $.

\item[\emph{(iii)}] There exists a unique (in distribution) L\'evy process
$Y_{\rho}$ such
that

$\mu = \mathcal{L}\Big(\int_{(0,1)}tdY_{\rho}(t)\Big)$, where
$\mathcal{L}(Y_{\rho}(1))=\rho \in \mathcal{U}^{<m-1>}$.

Moreover, in (ii) we have
$\mu_c=\mathcal{L}\Big(\int_{[c,1)}tdY_{\rho}(t)\Big)$, for
$0<c<1$.
\end{itemize}
\end{prop}
\emph{Proof.} For $m=1$, the above is just the Proposition 1. Now
suppose that the proposition is proved for $m$. If $\mu \in
\mathcal{U}^{<m+1>}$ then, by the definition (formula (1)),
$\rho_k \in \mathcal{U}^{<m>}$, for $k=1,2,...$ . For given
$0<c<1$, let us choose natural numbers $m_n$ such that $1\le
m_n\le n$ and $m_n/n\to c$, as $n \to \infty$. From (2) we have
\begin{equation}
\nu_n = T_{m_n/n}\nu_{m_n}^{\ast m_n/n}\ast
T_{1/n}(\rho_{m_{n}+1}\ast ...\ast\rho_n)^{\ast 1/n}.
\end{equation}
By Theorems 1.2 and 2.1 in Parthasarathy (1967), the second convolution
factor in (5) converges, say to $\mu_c$, which must be in
$\mathcal{U}^{<m>}$ by Corollary 1. Thus we get the factorization
(ii) for $m+1$, i.e., (i) implies (ii).

If (ii) holds we have a family $\mathcal{C}:=\{\mu_c : 0\le c\le
1\}\subset \mathcal{U}^{<m>}$, where $\mu_1=\delta_0$ and
$\mu_0=\mu$, from which we construct sequence $(\rho_k)$ as
follows
\[
\rho_1:=\mu \ \  \mbox{and} \   \rho_k:=T_k \mu_{(k-1)/k}^{\ast k}
\ \ \mbox{for} \ \ k \ge2.
\]
Using the factorization (ii) for $c=(k-1)/k$, then applying to
both sides the dilation $T_k$ and then raising to the (convolution)
power $k$, gives the equality $T_k \mu^{\ast k}=T_{k-1}\mu^{\ast
(k-1)}\ast \rho_k$, or in terms of Fourier transforms
\[
\hat\rho_{k}(y)= [\hat \mu(ky)]^k/[\hat \mu((k-1)y)]^{k-1}, \ \
\mbox{for} \ \  k \ge 2.
\]
Hence
\[
\rho_1 \ast \rho_2 \ast ... \ast \rho_n= T_n \mu^{\ast n} \ \
\mbox{i.e.,} \ \  \mu \in \mathcal{U}^{m+1},
\]
which completes
the proof that (ii) implies (i).

Since we have \begin{multline*} \mu =
\mathcal{L}(\int_{(0,1)}tdY_{\rho}(t))=
\mathcal{L}(\int_{(0,c)}tdY_{\rho}(t))\ast
\mathcal{L}(\int_{[c,1)}tdY_{\rho}(t))\\
= T_c(\mathcal{L}(\int_{(0,1)}tdY_{\rho}(t)))^{\ast c}\ast
\mathcal{L}(\int_{[c,1)}tdY_{\rho}(t)),
\end{multline*}
we infer that (iii) implies (ii). To prove the converse that (ii)
implies (iii) we proceed as in Jurek (1988), page 482 (formula
(3.1)) till page 484, taking $\beta=1$ and $Q=I$ (identity
operator). Thus we construct  process $Z(t)$ with independent
increments and cadlag paths such that $\mathcal{L}(Z(t))=
\mu_{e^{-t}}\in \mathcal{U}^{<m>}$. Because of Corollary 1 we
conclude that
\[
\tilde{Y}(t):=\int_{(0,t]}sdZ(s), \ \ \mbox{for} \ \ t \ge0,
\]
has increments with probability distributions in
$\mathcal{U}^{<m>}$.  All in all we have proved (iii).

\medskip
\begin{cor}
(a) The classes \ $\mathcal{U}^{<m>}, \ m=1,2,...$, of
the m-times s-selfdecomposable probability measures are closed convolution
subsemigropus , closed under convolution powers and the dilations
$T_d$.
\begin{multline}
(b) \ \  \mathcal{U}^{<m>}=
\mathcal{J}(\mathcal{J}(...(\mathcal{J}(ID)))), (\mbox{m-times
composition}), \\ L_{m+1} \subset \mathcal{U}^{<m+1>} \subset
\mathcal{U}^{<m>} \subset ID  , \ \ \mbox{for} \ \ m=0,1,2,... \ , \ \ \ \ \ \
\end{multline}
where $L_k$, $k=1,2,... $, are the convolution semigroups of
k-times selfdecomposable probability distributions.
\end{cor}
\emph{Proof.} Part (a) follows from the characterization (ii) in Proposition 2.
To prove that $\mathcal{U}^{<m>}$ are
closed we use Theorem 1.7.1 in Jurek and Mason (1993) or cf.
Chapter 2 in Parthasarathy (1967).

Part(b). Since $\mathcal{U} \subset ID$, therefore applying successively the random integral
mapping $\mathcal{J}$ to both sides gives the inclusion
$\mathcal{U}^{<m+1>}\subset \mathcal{U}^{<m>}$.
\newline
For the second inclusion $L_k \subset \mathcal{U}^{<k>}$, note that it
is true for $k=1$, cf. Corollary 4.1 in Jurek (1985). Assume it is true
$n$, i.e., $L_n \subset \mathcal{U}^{<n>}$ and let $\mu \in L_{n+1} \subset L_n$. Then for
any $0<c<1$ there exits $\nu_c \in L_n$ such that
\[
\mu = T_c \mu \ast \nu_c = T_c\mu^{\ast c}\ast \rho_c \ \mbox{with} \ \
\rho_c:=T_c\mu^{\ast (1-c)}\ast \nu_c \in L_n \subset \mathcal{U}^{<n>
},
\]
because, by the induction assumption, $\nu_c$ and $\mu$ are in
$L_n$. Consequently, by (ii), in Proposition 2, $\mu \in \mathcal{U}^{<n+1>}$
and this completes the proof.

\medskip
Our next aim is to describe m-times s-selfdecomposability in terms
parameters of infinitely divisible laws. Recall that each $ID$
distribution $\mu$ is uniquely determined by a triple: a shift
vector $a\in E$, a Gaussian covariance operator $R$, and a L\'evy
spectral measure $M$; we will write $\rho=[a,R,M]$. These are the
parameters in the L\'evy-Khintchine representation of the
characteristic function $\hat{\mu}$, namely
\begin{multline}
\rho \in ID$ iff $\hat{\rho}(y)= \exp(\Phi(y)), \ \ \mbox{where} \\
\Phi(y) : = i<y,a>+ 1/2<Ry,y>+  \\
\quad \int_{E \setminus\{0\}} [e^{i<y,x>}-1-i<y,x>1_{||x||\leq
1}(x)] M(dx),\ \ y \in E^{\prime};
\end{multline}
$\Phi$ is called the \emph{L\'evy exponent} of $\hat{\rho}$ (cf.
Araujo and Gin\'e (1980), Section 3.6). Furthermore, by the
\emph{L\'evy spectral function} of $\rho$ we mean the function
\[
L_{M}(D,r): = - M(\{x \in E : \|x\|>r \ \text{and} \ x\,||x||^{-1}
\in D \}),
\]
where $D$ is a Borel subset of unit sphere $S:=\{x: ||x||=1\}$ and
$r>0$. Note that $L_M$ uniquely determines $M$.

Since the L\'evy processes have infinitely divisible increments
(from the class $ID$) and $ID$ is a topologically closed
convolution semigroup, and also closed under dilations $T_a$ (a
multiplication of random variable by a scalar $a$), therefore the
random integrals $\int_{(a,b]}h(t)\,dY(r(t))$ have probability
distributions in $ID$ as well . If $[a_{h,r},R_{h,r},M_{h,r}]$
denotes the triple corresponding to the probability distribution
of the integral in question, and $[a,R,M]$ denotes the one
corresponding to the law of $Y(1)$ then (4) and (7) give the
following equation:
\begin{equation}
R_{h,r}=\big(\int_{(a,b]}h^2(t)dr(t)\big)\,R ,
\end{equation}
\begin{equation}
M_{h,r}(A)=\int_{(a,b]}M((h(t))^{-1}A)dr(t)\ \ \mbox{for} \ \ A\in
\mathcal{B}_0,
\end{equation}
and finally for the shift vector we have
\begin{multline}
a_{h,r}=\big(\int_{(a,b]}h(t)dr(t)\big)\,a \\
+\int_{E\setminus\{0\}}x\int_{(a,b]}h(t)[1_{B}(h(t)x)-1_B(x)]dr(t)M(dx).
\end{multline}
Specializing the above for the functions $h(t)=r(t)=t, \ 0 \le t
\le1$, and, for the simplicity of notations, putting
$\rho=[a,R,M]$ and $[a',R',M']= \mathcal{J}(\rho)$, we get from
(8)-(10) the following relations
\begin{equation}
R'= 1/3 R,
\end{equation}
\begin{equation}
M'(A)=\int_{(0,1)}M(t^{-1}A)dt, \ \ \mbox{for} \ \ A\in
\mathcal{B}_0,
\end{equation}
\begin{equation}
a'= \frac{1}{2}a + \int_{(0,1)}t\int_{1<||x||\le t^{-1}}x
M(dx)dt=\frac{1}{2}[a +\int_{\{||x|| > 1\}}x\,||x||^{-2}M(dx)].
\end{equation}
In order to get the second equality in (13) one needs to observe
that \newline $1_{\{t||x||\le1 \}}(x)=1_{\{0<t\le||x||^{-1}\}}(t)$ or
to change the order of integration. Thus
\begin{align*}
\int_{(0,1)}t \int_{\{1<||x||\le t^{-1}\}}x
M(dx)dt=\int_{\{||x|| > 1\}}\int_{(0,1)}tx 1_{\{t||x||\le1\}}(x)dtM(dx) \\
= \int_{\{||x|| > 1\}}x\int_0^{||x||^{-1}}tdtM(dx)= 1/2
\int_{\{||x|| > 1\}}x\,||x||^{-2}M(dx).
\end{align*}

\medskip
Now we may characterize the m-times s-selfdecomposable
distributions in terms of the triples in their L\'evy-Khintchine
formula.
\begin{prop}
For $m=1,2,...$, let $\rho=[a,R,M]$ and \newline $[a^{<m>},
R^{<m>},M^{<m>}]= \mathcal{J}^{m}(\rho)$ be m-times
s-selfdecomposable probability measures. Then
\begin{equation}
R^{<m>}= (1/3)^{m}R,
\end{equation}
\begin{equation}
M^{<m>}(A)=((m-1) !)^{-1}\int_{(0,1)}M(t^{-1}A)(- \ln t)^{m-1}dt, \
\ \mbox{for} \ \ A\in \mathcal{B}_0,
\end{equation}
\begin{equation}
a^{<m>}= (1/2)^{m}[a +\int_{\{||x|| >
1\}}x\,||x||^{-2}\,\sum_{j=0}^{m-1}\frac{(2\ln||x||)^j}{j!}\,M(dx).
\end{equation}
\end{prop}
\emph{Proof.} For $m=1$ the above are just the formulae (11)-(13).
\newline Assume that (14)-(16) holds form $m$. Since
$[a^{<m+1>}, R^{<m+1>},M^{<m+1>}]=\mathcal{J}([a^{<m>},
R^{<m>},M^{<m>}])$, therefore, by (11)-(13), we get
\begin{equation}
R^{<m+1>}= (1/3) R^{<m>},
\end{equation}
\begin{equation}
M^{<m+1>}(A)=\int_{(0,1)}M^{<m>}(t^{-1}A)dt, \ \ \mbox{for} \ \
A\in \mathcal{B}_0,
\end{equation}
\begin{equation}
a^{<m+1>}=\frac{1}{2}\,[a^{<m>} +\int_{\{||x|| >
1\}}x\,||x||^{-2}M^{<m>}(dx)\,].
\end{equation}
Obviously, by the induction assumption, we infer that (14) holds
for $m+1$. Similarly from (8) and (11), and from the change of the
order of integration, we get
\begin{align*}
M^{<m+1>}(A)= ((m-1)
!)^{-1}\int_{(0,1)}\int_{(0,1)}M(t^{-1}s^{-1}A)(-\ln
 t)^{m-1}dtds \\
=((m-1) !)^{-1}\int_{(0,1)}\int_{(0,t)}M(u^{-1}A)\frac{du}{t}(-\ln
t)^{m-1}dt
\\
 = ((m-1) !)^{-1}\int_{(0,1)}M(u^{-1}A)[\int_{(u,1)}\frac{1}{t}(-\ln
  t)^{m-1}dt]du \\
 =(m!)^{-1}\int_{(0,1)}M(u^{-1}A)(-\ln u)^{m}du,
\end{align*}
which proves (15).

In order to prove the formula for the shift, first note that by
(11) and by change of order of integration, we have
\begin{multline*}
w_m:=\int_{\{||x|| > 1\}}x\,||x||^{-2}M^{<m>}(dx) \\ = ((m-1)
!)^{-1}\int_{(0,1)}\int_E 1_{\{\{||x|| >
1\}\}}(tz)z\,||z||^{-2}\,t^{-1}(- \ln t)^{m-1}dt\,M(x)\\
= ((m-1)
!)^{-1}\int_{\{||z||>1\}}z\,||z||^{-2}\big[\int_{||z||^{-1}}^1
t^{-1}(-\ln t)^{m-1}dt\big]M(dz)\\ =
(m!)^{-1}\int_{\{||z||>1\}}z\,||z||^{-2}(\ln ||z||)^{m} M(dz),
\end{multline*}
for $m=1,2,...$ . Note that for $m=0$ the above formula gives the
second summand in (13). In terms of $w_m$, (19) gives the
recurrence relation
\[
a^{<m>}= 1/2(a^{<m-1>}+ w_{m-1}), \ \ \mbox{for} \ \ m=1,2,... ,
\]
where $a^{<0>}:= a$. Thus, if the formula for the shifts (16) holds
for $m$, then the above gives that it also holds for $m+1$, which
completes the proof the proposition.
\medskip

Let us recall that the functions
\begin{equation}
\Gamma(\alpha,x)=\int_x^{\infty}e^{-t}t^{\alpha -1}dt, \  \ x>0, \
\alpha >0,
\end{equation}
are called \emph{the incomplete gamma functions}. Simple
calculations shows that
\begin{equation}
\Gamma(m,x)= (m-1)!\,e^{-x}\,\sum_{j=0}^{m-1}\frac{x^j}{j!}, \ \
\mbox{for} \ \ m=1,2,... \  \ .
\end{equation}
Consequently, the formula (16) may be written as
\begin{equation}
a^{<m>}= (1/2)^m \,[\,a+ \frac{1}{\Gamma(m)}\int_{\{||x|| > 1\}}
x\, \Gamma(m, 2 \ln ||x||)\,M(dx)\,].
\end{equation}

Let us introduce $\emph{rescales of time}$ in the interval $(0,1)$
as follows
\begin{equation}
\tau_{\alpha}(t):= \frac{1}{\Gamma{(\alpha)}}\int_0^t (- \ln
u)^{\alpha-1} du,\ \ \ 0<t \le 1.
\end{equation}
Note that $\tau_{\alpha}$ is the cumulative probability
distribution function of the random variable
$g_{\alpha}:=e^{-G_{\alpha}}$, where $G_{\alpha}$ is the gamma
random variable with the probability density
$(\Gamma({\alpha}))^{-1}\,x^{\alpha -1}e^{-x}$ for $x>0$, and zero
elsewhere. Hence
\begin{multline}
\tau_{\alpha}(t)= P\{g_{\alpha} \le t\},\ \ \mbox{and} \ \
\mathbf{E}[g_{\alpha}^s]= (\Gamma(\alpha))^{-1}\int_0^1t^s(-\ln
t)^{\alpha-1}dt =(s+1)^{-\alpha}; \\
\int_0^c t^s\,d\tau_{\alpha}(t)= (s+1)^{-\alpha}\,
\frac{\Gamma(\alpha, -(s+1)\ln c)}{\Gamma(\alpha)}, \ \ \mbox{for}
\ \ s>0, \ 0<c<1,
\end{multline}
and (15) can rewritten as
\begin{equation}
M^{<m>}(A)=\int_{(0,1)}M(t^{-1}A)d\tau_m(t)=
\mathbf{E}[M(g_m^{-1}A)].
\end{equation}

\medskip
Now we can establish \emph{the random integral representation}
for the subclasses $\mathcal{U}^{<m>}$ of s-selfdecomposable
probability measures.
\begin{prop}
(a) The class $\mathcal{U}^{<m>}$ of m-times s-selfdecomposable
probability measures coincides with the class of probability
distributions of random integrals $\int_{(0,1)}t dY(\tau_m(t))$,
where $Y(\cdot)$ is an arbitrary L\'evy process.

(b) The class of Fourier transforms of measures from
$\mathcal{U}^{<m>}$ coincides with the class of functions
$\exp\mathbf{E}[\Psi(g_m y)], y\in E'$, where $\Psi$ is an
arbitrary L\'evy exponent of an infinitely divisible probability
measure and the random variable $g_m:=\exp (-G_m)$, with $G_m$
being the standard gamma random variable. In fact, the exponent
$\Psi$ is that of the random variable $Y(1)$ from (a).
\end{prop}
\emph{Proof.}
Let $[b_m,S_m,N_m]$ and $[a,R,M]$ are the triples
describing the probability distribution of the integral $\int_{(0,1)}t dY(\tau_m(t))$
and $Y(1)$, respectively. Then from (8)-(10) and (24) we have
\begin{multline*}
S_m =\big(((m-1) !)^{-1}\int_{(0,1)}t^2(-\ln t)^{m-1}dt \big)\,R=\frac{1}{3^m}\,R=R^{<m>}, \\
N_m(A)=((m-1) !)^{-1}\int_{(0,1)}M(t^{-1}A)(-\ln t)^{m-1}\,dt=M^{<m>}(A), \\
b_m=\big( ((m-1) !)^{-1}\int_{(0,1)} t (-\ln t)^{m-1}dt\big)a \, \\
+((m-1) !)^{-1}\int_{(0,1)}t \int_{1<||x||\le t^{-1}}x\,M(dx)(-\ln
t)^{m-1}\,dt
\end{multline*}
\begin{multline*}
=2^{-m}\,a +\int_{\{||x|| > 1\}}x
\big[((m-1)!)^{-1}\int_{(0,||x||^{-1})}t(-\ln
t)^{m-1}dt\big]\,M(dx)\\
=2^{-m}a + 2^{-m}\int_{\{||x|| > 1\}}x\,
[((m-1)!)^{-1}\,\Gamma(m,2\ln ||x||)\,]\, M(dx)= a^{<m>},
\end{multline*}
which completes the proof of the part (a).

For the part (b) we need to combine the formulae (4) (for $h(t)=t,
r(t)=\tau_m(t)$) and (14)-(16), and use (24).
\begin{cor}
A function $\phi : E'\to \Cset$ is a  Fourier transform of an
m-times s-selfdecomposable probability measure if and only if
there exist unique shift $a\in E$, a Gaussian covariance operator
$R$ and a L\'evy spectral measure $M$ such that
\begin{multline*}
\phi(y)= \exp \{i<y,a>+2^{-1}<Ry,y>+ \\
\int_{E\setminus\{0\}}[\,
\widehat{\big(\mathcal{L}(g_m)\big)}(<y,x>)-1-2^{-m}i<y,x>1_{B}(x)\,]M(dx)\},
\end{multline*}
where $g_m=e^{-G_{m}}$ and $G_m$ is the gamma random variable.
\end{cor}
\emph{Proof.} Use Proposition 4 together with the formula (4). Note that
there are not restrictions on a shift vector and a Gaussian covariance operator $R$.
Finally, for m=1 this is Theorem 2.9 in Jurek (1985).
\begin{rem}
Using the series representation of the exponential function and
(24) we get
\[
\widehat{\mathcal{L}(g_{\alpha}})(t)=\sum_{n=0}^{\infty}\frac{(it)^n}{n!
(n+1)^{\alpha}}, \ \ t\in \Rset.
\]
\end{rem}
Since the previous characterization, of m-times
s-selfdecomposability, has only a restriction  on the L\'evy
spectral measure, therefore we have a characterization of
$\mathcal{U}^{<m>}$ in terms of L\'evy spectral functions.
\begin{cor}
An infinitely divisible $[a,R,M]$ probability measure is m-times
s-selfdecomposable if and only if there exists a unique L\'evy
spectral measure $G$ such that
\[
L_M(D,r)= ((m-1)!)^{-1}r\int_r^{\infty}(\ln w- \ln
r)^{m-1}L_G(D,w)\frac{dw}{w^{2}}, \] for all sets $D$ and  all
$r>0$, or equivalently
\[
L_{M}(D,r)=r\int_{r}^{\infty}x_{m-1}^{-1}\int_{x_{m-1}}^{\infty}
... \,
x_1^{-1}\int_{x_1}^{\infty}w^{-2}L_G(D,w)dw\,dx_1\,...\,dx_{m-1},
\]
for all sets $D$ and all $r>0$.
\end{cor}
\emph{Proof.} In view of the Proposition 3 we have that
$M=G^{<m>}$ for a unique L\'evy spectral measure $G$, and (15)
gives the first part of the corollary. Since the relation (18), in
terms L\'evy spectral functions, reads
\[
L_{M^{<j>}}(D,r)=
r\int_r^{\infty}L_{M^{<j-1>}}(D,x_{j-1})\frac{dx_{j-1}}{x_{j-1}^2},
\ \ r>0, \ \mbox{with} \ \ M^{<0>}:= G,
\]
for $j=1,2...$, therefore the inductive argument proves the second
part of the corollary.
\begin{cor}
In order that a L\'evy spectral measure $G$ to be a L\'evy spectral
measure of an m-times s-selfdecomposable probability measure, it
is necessary and sufficient that its L\'evy spectral functions
$r\to L_G(D,r)$ are m-times differentiable, except at countable
many points $r$, and the function
\[
L(D,r)=(\mathcal{A}^m(L_G(D,\cdot)))(r)\ \ \mbox{is a L\'evy
spectral function. }
\]
The operator $\mathcal{A}^{m}$ is the m-time composition of the
linear differential operator $\mathcal{A}$, which is  defined as follows
\[
   \ \ \ \ (\mathcal{A}(h))(x):= xh'(x)-h(x),
\]
for once differentiable real-valued functions $h$ defined on
$(0,\infty)$.
\end{cor}
\emph{Proof.} If measures $M'$ and $M$ are related as in (12) then
their corresponding spectral functions (tails)  $L_{M'}$ and $L_M$
satisfy equality
\[
L_{M'}(D,r)=\int_{(0,1)}L_M(D,r/t)dt=r\int_r^\infty
w^{-2}L_M(D,w)dw.
\]
Hence $L_{M'}$ is at least once differentiable (except on a
countable set) and
\begin{equation}
-L_M(D,r)= r \frac{d}{dr}L_{M'}(D,r) -
L_{M'}(D,r)=(\mathcal{A}(L_{M'}(D,\cdot))(r).
\end{equation}
Moreover, by Theorem 1.3 in Jurek (1985), we have the left-hand
side is a L\'evy spectral measure (on Banach space) if and only if
so is a measure on the right-hand side. Because of the recurrence
equation (18) we have proved the corollary.
\begin{rem}
The random integral mapping $\mathcal{J}$ is defined on infinitely
divisible measures $\rho=[a,R,M]$. If one assume that the formula
(12) defines the mapping $\mathcal{J}$ on the measure $M$ or its spectral
function $L_M$, then $\mathcal{A}$ may be viewed as its inverse
mapping.
\end{rem}
Before the next characterization, of the class $\mathcal{U}^{<m>}$
distributions, let us recall that a L\'evy exponent is just the
logarithm of an infinitely  divisible Fourier transform; cf.
formula (7). Let us note that, if $\Psi$ is the L\'evy exponent of
$\rho$ and $\Phi$ is that of $\mathcal{J}(\rho)$, then (3) and (4)
give the following
\[
\Phi(ty)=\int_{(0,1)}\Psi(sty)ds=
\frac{1}{t}\int_{(0,t)}\Psi(sy)ds,
\]
and consequently
\[
\Psi(y)= \Phi(y)+ d(\Phi(ty))/dt|_{t=1}.
\]
With these equalities and the recursive relation between classes
$\mathcal{U}^{<m>}$ we have
\begin{cor}
A function $\Phi:E'\to\Cset$ is a L\'evy exponent of an
m-times s-selfdecomposable probability measure if and only if
there exists a unique L\'evy exponent $\Psi$ such that the
function
\[
E' \ni y \to \mathcal{D}^m(\Psi)(y) \ \ \mbox{ is a L\'evy
exponent.}
\]
The operator $\mathcal{D}^m$ is the m-time composition of the
following linear differential operator
\[
(\mathcal{D}g)(y):= g(y)+ d(g(ty))/dt|_{t=1},
\]
where $ g:E' \to \Cset$ is once differentiable in each direction
$y\in E'$ and $t\in \Rset$.
\end{cor}
\medskip
Note that in a particular case one has $(\mathcal{D}g)(y):= g(y)+
y\,dg(y))/dy$, when $y \in E'=\Rset$ and it differs from
$\mathcal{A}$ in Corollary 5 only by a sign.
\begin{rem}
If ones defines $\mathcal{J}$ on L\'evy exponents by (4) then the operator
$\mathcal{D}$ can be viewed as its inverse, i.e., $\mathcal{D}=\mathcal{J}^{-1}$,
on L\'evy exponents on a Banach space.
\end{rem}
\begin{prop}
A probability measure $\mu=[a,R,M]$ is completely
s-selfdecomposable, i.e., $\mu \in \mathcal{U}^{<\infty>}:=
\bigcap_{m=1}^{\infty} \mathcal{U}^{<m>}$ if and only if there
exists a unique bi-measure $\sigma(\cdot,\cdot)$ on $S\times(0,2)$
such that
\begin{equation}
M(A \cdot D)=\int_{(0,2)}\int_D \int_A \,w^{-(z+1)}\,dw\,\sigma(du,dz)=\int_{(0,2)}\int_A
\,w^{-(z+1)}\,dw\,\sigma(D,dz),
\end{equation}
where $A\cdot D:= \{x \in E: x/||x|| \in D , \ \ ||x|| \in A \}$
and for each Borel $D\subset S$, $\sigma(D,\cdot)$ is a finite
Borel measure on the interval $(0,2)$ and for each Borel subset
$A\subset (\epsilon,\infty)$ for some $\epsilon >0$, $\sigma(\cdot, A)$ is a finite Borel
measure on the unit sphere $S$. Moreover, we have that
\[
\int_{(0,2)}\int_S |<y,u>|^2 \frac{1}{2-z}\,\sigma(du,dz) <\infty,
\]
for all $y\in E'$.
\end{prop}
\emph{Proof.} If $\mu=[a,R,M]$ is completely s-selfdecomposable then by Proposition 3 or Corollary 4,
for each $m$ there exists a unique L\'evy measure $G$ such that
\[
M(A) =((m-1)!)^{-1}\int_{(0,1)}G(t^{-1}A)(- \ln t)^{m-1}dt, \ \mbox{for} \ \ A\in \mathcal{B}_0.
\]
or for all $D$ and $r>0$
\[
L_{G}(D,r)=r\int_{r}^{\infty}x_{m-1}^{-1}\int_{x_{m-1}}^{\infty}
... \,
x_1^{-1}\int_{x_1}^{\infty}w^{-2}L_G(D,w)dw\,dx_1\,...\,dx_{m-1}.
\]
Hence, for each set $D$, the functions $f(x):= - e^{-x}L_G(D,e^x),\ \ x \in \Rset$, are m-times differentiable
and
\[
(-1)^m d^m f(x)/dx^m = - e^{-x}\,L_G(D,e^x)\ge 0.
\]
In other words, $f$ is completely monotone and by Bernstein's Theorem, there exists a unique finite Borel measure
$\sigma^{\sim}(D,\cdot)$ on $(0,\infty)$ such that
\begin{equation}
f(x)=\int_0^{\infty}e^{-xz}\sigma^{\sim}(D,dz),\ \ \mbox{i.e.,} \ \ L_M(D,r)=
 -\int_0^{\infty}\frac{1}{r^{z-1}}\sigma^{\sim}(D,dz).
\end{equation}
Since L\'evy spectral functions vanish at $\infty$ and
$\sigma^{\sim}(D,\cdot)$ are finite measures, therefore they must be
concentrated on half-line $(1,\infty)$. Consequently, from (28) we
get
\begin{multline*}
M([r,s) \cdot D))=L_M(D,r)-L_M(D,s) =\int_1^{\infty}(z-1)\int_r^s
\frac{1}{w^z}\,dw\, \sigma^{\sim}(D,dz),
\end{multline*}
for all
$0<r<s<\infty$ and all Borel sets $D\subset S$. Since for $y\in
E'$, $(\pi_y\,M)(C):=M(\{x \in E : <y,x>\in C\})$, for Borel
subsets $C$ in $\Rset$, are L\'evy measure on real line therefore
\begin{multline*}
\int_{\{||x||\le 1\}}|<y,x>|^2\,M(dx) \le  \int_{\{x\in E:\,|<y,x>|\le1 \}}|<y,x>|^2 M(dx) \\
=  \int_{\{|t|\le1\}}t^2\,(\pi_y\,M)(dt) <\infty \,\, .
\end{multline*}
On the other hand, using (28) the integral
\begin{multline*}
\int_{\{0<||x||\le1\}}|<y,x>|^2\,M(dx)=\int_{(0,1]\cdot S}|<y,u>|^2\,t^2\,M(du \cdot dt) \\
= \int_S\int_0^1|<y,u>|^2\,t^2\int_1^{\infty}\frac{z-1}{t^z}\,\sigma^{\sim}(du,dz)\,dt\\
=\int_S\int_1^{\infty}|<y,u>|^2(z-1)\,[\int_0^1
t^{-(z-2)}dt\,]\,\sigma^{\sim}(du,dz),
\end{multline*}
is finite only if $z<3$, because $\sigma^{\sim}(D,\cdot)$ are
finite measures. Changing the variable and putting $\sigma(D,dz):=
z\sigma^{\sim}(D,dz+1),\ \ 0<z<2,$ we obtain the formula (27)
together with the integrability condition. Thus the necessity part
of the proposition is proved.

For the converse, let $\rho=[a,R,M]$ with the spectral measure $M$
of the form in (27). Hence, by (12), $\mathcal{J}(\rho)$ has
spectral measure
\[
M'(A\cdot D)=\int_{(0,1)} \int_{(0,2)}\int_{s^{-1}A}
\frac{1}{w^{z+1}} \,dw\,\sigma(D,dz)ds=\int_{(0,2)}\int_A
\frac{1}{x^{z+1}}dx\,\sigma_1(D,dz),
\]
where $\sigma_1(D,dz):= (z+1)^{-1}\sigma(D,dz)$ is another finite
measure on the interval $(0,2)$ and therefore $\mathcal{J}(\rho)$
has the L\'evy spectral measure of the form (27) again.
Consequently, $\rho \in \mathcal{U}^{<m>}$ for all $m$, and thus
the sufficiency of (27) is completed.

\medskip
Let put
\[
\Sigma:=\{\rho=[a,R,M_{\sigma}]: M_{\sigma} \ \ \mbox{\emph{is of
the form (27)}} \},
\]
i.,e., $\Sigma$ is the more explicit description of
$\mathcal{U}^{<\infty>}$. Further, let us recall that
\[
\int_E \log(1+||x||)\,\rho(dx)\ \ < \infty \ \ \mbox{iff} \ \
\int_{\{||x||>1\}}\log ||x||\,M_{\sigma}(dx)<\infty,
\]
By the formula (27), the last integral is equal
\begin{multline*}
- \int_1^{\infty} \log t\,dL_{M_{\sigma}}(S,t)=
\int_1^{\infty}L_{M_{\sigma}}(S,t)\,t^{-1}\,dt \\
= \int_1^{\infty}\int_{(0,2)}t^{-z-1}\,\sigma(S,dz)\,dt=
\int_{(0,2)}z^{-1}\,\sigma(S,dz) \ <\infty.
\end{multline*}
Similar integrability formulas hold for functions
$g_k(x):=\log^k(1+||x||)$ and L\'evy measures $M$. Recall that the
integrability condition of $g_k$  appears in the random integral
representation for the class $L_k$.
\begin{cor}
The class of completely s-selfdecomposable probability measures
coincides with the class of completely selfdecomposable ones,
i.e., $\mathcal{U}^{<\infty>}=L_{\infty}$.
\end{cor}
\emph{Proof.} If $\rho=[a,R,M]$ and
$\mathcal{I}(\rho):=\mathcal{L}\big(\int_0^{\infty}e^{-t}\,dY_{\rho}(t)\big)=
[A^0,R^0,M^0]$, then
\[
\int_{||x||>1} \log ||x|| M(dx)<\infty \ \ \mbox{and} \ \
M^{0}(A):=\int_0^{\infty}M(e^sA)ds,
\]
for all sets $A \in \mathcal{B}_0\,$ cf. Jurek (1985), p.603 or
Jurek and Mason (1993), p.120. Simple calculation show that
\begin{multline*}
(M_{\sigma})^0(A \cdot D)=\int_0^{\infty}\int_{(0,2)}
\int_{e^tA}\frac{1}{w^{z+1}}dw\,\sigma(D,dz)dt \\
= \int_{(0,2)}\int_A v^{-(z+1)}\int_0^\infty
\,e^{-tz}dt\,dv\sigma(D,dz)= \int_{(0,2)}\int_A
\frac{1}{v^{z+1}}dv\,\sigma_2(D,dz),
\end{multline*}
where $\sigma_2(D,dz):= z^{-1}\sigma(D,dz)$ is another finite
measure on $(0,2)$ because of the logarithmic moment assumption.
This shows that $\mathcal{I}(\Sigma)\subset \Sigma$. Since $L_k=
\mathcal{I}(\mathcal(...(\mathcal{I}(ID_{log^k}))))$, (k-times
composition of $\mathcal{I}$ and $ID_{log^k}$ denotes the class of
infinitely divisible measures with finite $\log^k$-moments) we
infer that $\Sigma\subset L_k$, for $k=1,2... \,$. Consequently,
$\Sigma \subset L_{\infty}:= \cap_{k=1}^{\infty}L_k \subset
\mathcal{U}^{<\infty>}=\Sigma$, which completes the proof.
\begin{rem}
Measures $M_{\sigma}$ are mixtures of L\'evy measures of stable
laws. The mixture is done with respect to the exponents $p\in
(0,2)$. Since Fourier transforms of p-stable measures are known
explicitly we can have analogous formulas for completely
s-selfdecomposable measures; cf. a similar result (on the real
line) for $L_{\infty}$ in Urbanik (1973), or Thu (1986), or Sato
(1980) or Jurek (1983),Theorem 7.2.
\end{rem}
\medskip
\textbf{3. Concluding remarks and two examples.}
\newline
\textbf{A).} The classes $\mathcal{U}^{<m>}$ were
introduced by an inductive procedure and thus we have the natural
index $m$. For a positive non-integer $\alpha$ one may proceeds as
in Thu (1986) using the fractional calculus. However, we may
utilize our random integral approach and define
\begin{equation}
\mathcal{U}^{<\alpha>
}=\{\mathcal{L}\big(\int_{(0,1)}t\,dY_{\rho}(\tau_{\alpha}(t))\big):
\rho\in ID \},
\end{equation}
where $Y_{\rho}(\cdot)$ is a L\'evy process with
$\mathcal{L}(Y_{\rho}(1))=\rho$. Equivalently, we have
\[
[a^{<\alpha>},R^{<\alpha>},M^{<\alpha>}]=\mathcal{J}^{\alpha}(\rho)=
\mathcal{L}\big(\int_{(0,1)}t\,dY_{\rho}(\tau_{\alpha}(t))\big),
\]
where
\begin{multline}
R^{<\alpha>}= 3^{-\alpha}R , \ \ \
M^{<\alpha>}(A)=\int_{(0,1)}M(t^{-1}A)d\tau_{\alpha}(t), \ \ A\in
\mathcal{B}_{0}, \\
a^{<\alpha>}= 2^{-\alpha}\,[\,a+
\frac{1}{\Gamma(\alpha)}\int_{\{||x|| > 1\} }x\,\Gamma(\alpha,
2\ln ||x||)\,M(dx)],
\end{multline}
cf. (14), (15) and for the shift vector (16) with (21),(22) and
(24).

Furthermore, for any continuous and bounded $f$ on $(0,\infty)$
and gamma random variables  $G_{\alpha}$ and $G_{\beta}$ we have
\[
\int_0^{\infty}\int_0^{\infty}f(x+y)\mathcal{L}(G_{\alpha})(dx)\mathcal{L}(G_{\beta})(dy)=
\int_0^\infty f(z)\mathcal{L}(G_{\alpha+\beta})(dz),
\]
i.e., gamma distributions form an one-parameter convolution
semigroup of measures on $(0,\infty)$ with the addition .
Consequently, for any continuous bounded $h$ on interval $(0,1)$
\begin{equation*}
\int_{(0,1)}\int_{(0,1)}h(st)d\tau_{\alpha}(t)d\tau_{\beta}(s)=
\int_{(0,1)}h(u)d\tau_{\alpha+\beta}(u),
\end{equation*}
thus $\tau_{\alpha}$ form one parameter semigroup of measures on
$(0,1)$ with the multiplication. Hence we infer that
\begin{cor}
For any positive $\alpha$ and $\beta$ we have
\begin{itemize}
\item[\emph{(a)}] $\mathcal{J}^{<\alpha +\beta>}=
\mathcal{J}^{\alpha}(\mathcal{J}^{\beta})$,

\item[\emph{(b)}] if $\alpha <\beta$ then
$\mathcal{U}^{<\beta>} \subset \mathcal{U}^{<\alpha>}$.
\end{itemize}
\end{cor}
\medskip

\textbf{B).} In this subsection we consider only $\Rset$-valued
random variables or Borel measures on the real line.
\newline
Because of the inclusion $L \subset \mathcal{U}$ each
selfdecomposable distribution is an example of s-selfdecomposable
one. On the other hand, by Proposition 3 in Iksanow, Jurek and
Schreiber (2002), selfdecomposable distributions of random
variables of the form $X:=\sum_{k=1}^{\infty}a_k \eta_k$, where
$\eta_k$'s are independent identically distributed Laplace (double
exponential) random variables and $\sum_k a_k^2<\infty$, have the
background driving probability measures $\nu \in \mathcal{U}$.
Furthermore, by Proposition 3 in Jurek (2001) we have that
\begin{equation}
\hat{\nu}(t)=\exp\,[t \phi'_{X}(t)/\phi_{X}(t)\,], \ \  \ \ t \in \Rset.
\end{equation}
In Jurek (1996) it was noticed that $\phi_S(t)=t/(\sinh t)$ ("S"
stands for the hyperbolic 'sine') and $\phi_{C}(t):= 1/(\cosh t)$
( "C" stands for the hyperbolic 'cosine') are the characteristic
functions of random variables of the above series form $X$. Using
(31) we conclude
\[
\psi_S(t):= \exp(1-t\coth t), \ \psi_C(t):=\exp(-t\tanh t) \ \
\mbox{are class $\mathcal{U}$ char. f.}
\]
Thus both are characteristic functions of integrals (3).
Furthermore from Corollary 6 we have that
\begin{multline}
\mathcal{D}(\log \psi_S(t))=1-2\coth t +t^2/(\sinh^2t) \ \
\mbox{and}\\
\mathcal{D}(\log \psi_C(t))=-2t\tanh t - t^2/(\cosh^2t) \ \
\mbox{are L\'evy exponents.}
\end{multline}
It might be worthy to mention here that $\phi_S(t)\cdot \psi_S(t)$
is a characteristic function of \emph{a conditional L\'evy's
random area integral}; cf. L\'evy (1951) or Yor (1992) and Jurek
(2001). Similarly, $(\phi_C(t)\cdot\psi_C(t))^{1/2}$ is a
characteristic function of \emph{an integral functional of
Brownian motion}; cf. Wenocur (1986) and Jurek (2001), p. 248.

Recently in Jurek and Yor (2002) the probability distributions
corresponding to both $\psi_S$ and $\psi_C$ were expressed in
terms of squared Bessel bridges. Also both functions viewed as the
Laplace transform in $t^2/2$ can be interpreted as the hitting
time of 1 by the Bessel process starting from zero; cf. Yor
(1997), p. 132. At present we are not aware of any stochastic
representation for the analytic expressions in (32). Finally, it
seems that the operators $\mathcal{A}^m$ may be related to some
Markov processes.

\medskip
\noindent Institute of Mathematics, University of Wroc\l aw,
50-384 Wroc\l aw, Poland. \newline [E-mail:
zjjurek@math.uni.wroc.pl ]\newline
[Currently: Department of Mathematics, Wayne State University, Detroit, MI 48202, USA.]


\medskip
\begin{thebibliography}{29}

\bibitem{AG} A. Araujo and E. Gin\'e, \emph{The Central Limit
Theorem for Real and Banach Valued Random Variables}, Wiley, New
York, 1980.

\bibitem{IJS}  A. M. Iksanov, Z. J. Jurek and B.M. Schreiber, A new
factorization property of the selfdecomposable probability
measures. Preprint 2002.

\bibitem{JJM} J. Jacod, A. Jakubowski and J. M\'emin, \emph{About
Asymptotic Errors in Discretization}, Prepublication No. 661,
Laboratoire de Probabilit\'es et Mod\`eles Al\'eatoires, Univ.
Paris VI, May, 2001.

\bibitem{J77} Z. J. Jurek, Limit distributions for sums of
shrunken random variables. In \emph{Second Vilinius Conf. Probab.
Theor. Math. Stat.}, Abstract of Communications \textbf{3}, 95-96,
1977.

\bibitem{J81} Z. J. Jurek, Limit distributions for sums of
shrunken random variables in Hilbert spaces, \emph{Dissertationes
Math.( Rozprawy Matematyczne)}, \textbf{185}, PWN Warszawa, 1981.

\bibitem{J83} Z. J. Jurek, Limit distributions and one-parameter
groups of linear operators on Banach spces, \emph{J. Multivariate
Anal.} \textbf{13}, 578--604, 1983.

\bibitem{J85} Z. J. Jurek, Relations between the $s$-selfdecomposable and
selfdecomposable measures, \emph{Ann. Probab.} \textbf{13},
592--608, 1985.

\bibitem{J96} Z. J. Jurek, Series of independent exponential random
variables. In \emph{Proc.7th Japan-Russia Symposium on Probab.
Theor. Math. Statist., Tokyo 1995}, S. Watanabe, M. Fukushima, Yu.
V. Prohorov and A. N. Shiryaev, eds., World Scientific, Singapore,
174--182, 1996.

\bibitem{J01} Z. J. Jurek, Remarks on the selfdecomposability and new
examples, \emph{Demonstratio Mathematica} \textbf{34}, 241--250,
2001.

\bibitem{JM} Z. J. Jurek, J.D. Mason, Operator-limit distributions
in probability theory, \emph{John Wiley \& Sons,Inc.}, New York,
1993.

\bibitem{JV} Z. J. Jurek, W. Vervaat, An integral representation
for selfdecomposable banach space valued random variables,
\emph{Z. Wahrscheinlichkeitstheor. Verw. Geb.} \textbf{62},
247-262, 1983.

\bibitem{JY} Z. J. Jurek, M. Yor, Selfdecomposable laws assoocitaed with
hyperbolic functions. Preprint 2002.

\bibitem{L} P. L\'evy, Wiener's random functions, and other
Laplacian random functions. In \emph{Proc. 2nd Berkeley Symposium
Math. Stat. Probab.}, Univ. California Press, Berkeley, 171--178,
1951.

\bibitem{O'C}  T. A. O'Connor, Infinitely divisible distributions with
unimodal L\'{e}vy spectral functions, \emph{Ann. Probab.}
\textbf{7}, 494--499, 1979.

\bibitem{P} K. R. Parthasarathy, Probability measures on metric
spaces, \emph{Academic Press}, New York, 1967.

\bibitem{S} K. Sato, Class $L$ of multivariate distributions and its
subclasses, \emph{J. Multivarite Anal.} \textbf{10}, 201--232,
1980.

\bibitem{T} N. Thu, An alternative approach to multiply
selfdecomposable probability meaures on Banach spaces,
\emph{Probab. Th. Rel. Fields} \textbf{72}, 35-54, 1986.

\bibitem{W} M. Wenocur, Brownian motion with quadratic killing
and some implications, \emph{J. Appl. Probab.} \textbf{23},
893--903, 1986.

\bibitem{Y} M. Yor, \emph{Some Aspects of Brownian Motion, Part I:
Some Special Functionals}, Birkh\"auser, Basel, 1992.

\bibitem{Y2} M. Yor, \emph{Some Aspects of Brownian Motion, Part
II: Some recent martingale problems,} Birkh\"auser, Basel,
1997.

\bibitem{U} K. Urbanik, \emph{Limit laws for sequences of normed
sums satysfying some stability conditions.} In: Multivariate
Analysis III, P.R. Krishnaiah Ed., 225-237, Academic Press, New
York, 1973.

\end{thebibliography}
\end{document}